\documentclass[11pt]{article}
\usepackage{amsfonts,enumerate,amsmath}

 \newenvironment{dedication}
        {\vspace{3ex}\begin{quotation}\begin{center}\begin{em}}
        {\par\end{em}\end{center}\end{quotation}}

\newtheorem{theorem}{Theorem}

\def\C{{\mathbb C}}
\def\R{{\mathbb R}}
\def\Z{{\mathbb Z}}
\def\Q{{\mathbb Q}}
\def\N{{\mathbb N}}
\def\M{{\mathcal M}}

\def\beq{\begin{equation}}
\def\eeq{\end{equation}}

\begin{document}

\title{The spaces of non-contractible closed curves in compact space forms}
\author{I.A. Taimanov
\thanks{Sobolev Institute of Mathematics, 630090 Novosibirsk, Russia, and Novosibirsk State University, 630090 Novosibirsk, Russia; e-mail: taimanov@math.nsc.ru.
The work was supported by RSF (grant 14-11-00441). }
}
\date{}
\maketitle
\begin{dedication}
\hfill{To the memory of D.V. Anosov}
\end{dedication}

\section{Introduction}

The study of the periodic problem for Finsler geodesics was initiated by Anosov in \cite{Anosov74}.
He explained his interest in the study of Finsler geometry by pointing out that ``it leads to a wide class of dynamical systems which admits
an application of geometrical notions and reasonings for formulating problems and for the study of them''.
For instance, in early 1980s the study of variational problems for magnetic geodesic flows,
which under certain conditions are particular cases of Finsler geodesic flows, had been started in \cite{Novikov} (see, also,  \cite{T0}). 

In \cite{Anosov74} Anosov claimed that, in difference with the Riemannian case
where, by the Lyusternik--Schnirelman theorem, there exists at least three nonselfintersecting closed geodesics on  the two-sphere,
for irreversible Finsler metrics ``one can guarantee only the existence of two closed geodesics''.
This fundamental result had been proved only recently by Bangert and Long \cite{BL} and the proof strongly relies on
the index iteration formulas derived by Long \cite{Long}. The estimate can not be improved due to the Katok example \cite{Katok}.

Recently Xiao and Long studied the topological structure of non-con\-trac\-tible loop spaces for odd-dimensional projective spaces
computing, in particular, the equivariant cohomology with $\Z_2$-coefficients of the path spaces \cite{XL}
and together with Duan applied these results to proving the existence of at least two geometrically distinct non-contractible
closed geodesics for irreversible bumpy Finsler metrics on $\R P^3$ \cite{DLX}.

In this article we demonstrate how to use a result from \cite{Ballmann,T86} for calculating the rational
equivariant cohomology of non-contractible loop spaces for the compact space forms. 
We also show how to use these calculations for  establishing the existence of closed geodesics.

\section{The path spaces}
\label{sec2}

Let $M^n$ be a closed Riemannian manifold.

Let us denote by
$\Lambda(M^n)=H^1(S^1,M)$ the space of $H^1$-maps
$$
\gamma: [0,1] \to M^n, \ \ \ f(0)=f(1),
$$
of a circle $S^1 = \R/\Z$ into $M^n$, by $\Omega_x (M^n)$ the subspace of
$\Lambda(M^n)$ formed by loops starting and ending at $\gamma(0)=\gamma(1) =x \in M^n$, and by
$\Pi^+(M^n)$ and $\Pi(M^n)$ the quotients of $\Lambda(M^n)$ with respect to the $SO(2) (=S^1)$-action:
$$
\varphi \cdot \gamma(t) = \gamma(t+\varphi), \ \ \ \varphi \in S^1 = \R/\Z,
$$
and the $O(2)$-action respectively. Here the $O(2)$ action is the extension of the $SO(2)$-action by
the involution
$$
\sigma \cdot f(t) = f(-t).
$$
The space $\Lambda(M)$ is a Hilbert manifold and the vibration  $H^1(S^1,TM) \to H^!(S^1,M) = \Lambda(M)$ is 
the tangent bundle to it. A detailed exposition of the topology and the atlas of Hilbert manifold on  $\Lambda (M)$ one can find in \cite{K,Anosov80}.
The space $\Lambda$ which was introduced into  the calculus of variations in
the middle of 1960s has many functorial properties.
In particular,

\begin{enumerate}
\item
if $f: M \to N$ is a smooth map, then the induced map
$\Lambda f: \Lambda(M) \to \Lambda(N)$ is a smooth map of Hilbert manifolds.
Moreover, if $f_s: M \times [a,b] \to N$ is a smooth homotopy, then $\Lambda f_s$ is also a smooth homotopy;

\item
the map $g: \Lambda M \to \Lambda M$ which corresponds to a curve the same geometrical curve
parameterised proportionally to the arc-length, is continuous. Therewith one may assume that the condition
$( g\cdot \gamma)(0) = \gamma(0)$ also holds;

\item
the aforementioned map $g$, which satisfies $(g \cdot \gamma)(0) = \gamma(0)$, is homotopical to the identity and moreover
such a homotopy can be chosen to be $SO(2)$- or $O(2)$-invariant.
\end{enumerate}

\noindent
The proofs of these statements are given in \cite{K,Anosov80}.

Let us denote by
$L(M)$ the space
$$
L(M) = g(\Lambda(M)),
$$
and by $P^+(M)$ and $P(M^n)$ the following quotient-spaces:
$$
P^+(M) = L(M)/SO(2), \ \ \ \ \  P(M)= L(M)/O(2).
$$
These spaces are deformation retracts of $\Lambda(M)$, $\Pi^+(M)$ and $\Pi(M)$ respectively and therefore are
homotopically equivalent to them.

Geometrically closed geodesics of the metric $g_{ik}\,dx^i\,dx^k$ are
the extremals of the energy functional
$$
E(\gamma) = \frac{1}{2} \int_\gamma |\dot{\gamma}|^2\, dt, \ \ \ E: \Lambda(M) \to \R,
$$
and of the length functional
$$
S(\gamma) = \int |\dot{\gamma}|\,dt, \ \ \ S: L(M) \to \R,
$$
where
$$
|\dot{\gamma}| = \sqrt{g_{ik}\dot{\gamma}^i\dot{\gamma}^k}.
$$
The Euler--Lagrange  equations for the energy functional imply that the parameter on an extremal
has to be proportional to the arc-length and to every closed geodesic there corresponds
a pair of $S^1$-families (an $O(2)$-orbit)  of extremals. For the length functional
every reparameterisation of an extremal is again an extremal.
Therefore we have to fix a parameter which is proportional to the arc-length to obtain again
a pair of $S^1$-families of extremals in $L(M)$.

A manifold $M$ with a function $F(x,\dot{x})$ defined on its tangent bundle is called a Finsler manifold if

\begin{enumerate}
\item
$F(x,\dot{x}) \geq 0$ and it vanishes if and only if $\dot{x}=0$;

\item
$F(x,\lambda\dot{x}) = \lambda F(x,\dot{x})$ for all $\lambda >0$;

\item
the unit spheres $\{F(x,\dot{x})=1\}$ are convex and their curvatures are positive with respect to
the Euclidean metrics in the tangent spaces $T_x(M)$.
\end{enumerate}

On such a manifold for every smooth path $\gamma(t), a \leq t \leq b$, there is defined its Finsler length
$$
S(\gamma) = \int_a^b F(\gamma,\dot{\gamma})\,dt.
$$
Since $F(x,\dot{x})$ is homogeneous of the first order in $\dot{x}$, every reparameterisation of
an extremal of $S$ is again its extremal. Therefore the variational problem for closed Finsler geodesics is
better to pose for the functional
$$
S(\gamma) = \int_\gamma F(\gamma,\dot{\gamma})\,dt, \ \ \ S: L(M) \to \R.
$$
In this case to every extremal there corresponds an $SO(2)$-orbit of extremals in $L(M)$.
If the Finsler metric is reversible, i.e. $F(x,\dot{x}) = F(x,-\dot{x})$, in particular, if $F(x,\dot{x}) = |\dot{x}|$ for some Riemannian metric,
then every not-one-point extremal generates an $O(2)$-orbit consisting of extremals.

To reduce the amount of critical points  corresponding geometrically to the same closed curve
we have to consider the length functionals $S$ on the spaces $P^+(M)$ (for irreversible metrics) and
$P(M) = P^+(M)/\sigma$ (for reversible metrics). It is clear that

{\sl to every extremal closed curve
there corresponds a unique critical point of $S$ in $P^+(M)$ (for irreversible metrics) and in
$P(M) = P^+(M)/\sigma$ (for reversible metrics).}

Two closed geodesics are called {\it distinct} if they are not both iterates $\gamma^n$ of the same closed curve $\gamma$ where  
$n \in \{1,2,\dots\}$ for irreversible metrics and $n \in \{\pm 1, \pm 2,\dots\}$ for reversible metrics.

\section{Rational homotopy of the path spaces}

Let $h \in \pi_1(M,x_0)$ and let $[h]$ be the corresponding free homotopy class of closed  curves: $[h] \in [S^1,M]$.
We denote by
$$
\Lambda M[h] \subset \Lambda M \ \ \ \mbox{and} \ \ \ LM[h] \subset LM
$$
the connected components of $\Lambda M$ and $LM$ consisting of curves from $[h]$.

Let
$h$ be realized by a map $\omega: [0,1] \to M$ with $\omega(0) = x_0$, and
let $h_i$ be the automorphism
$$
h_i: \pi_i(M,x_0) \to \pi_i(M,x_0)
$$
corresponding to the standard action of $h \in \pi_1$ on $\pi_i$.

The following theorem was proved independently in \cite{Ballmann} and \cite{T86}.

\vskip4mm

{\bf Theorem A} \ \cite{Ballmann,T86} \
{\sl The mapping
$$
\pi: \Lambda M \to M, \ \ \ \pi(\gamma) = \gamma(0),
$$
which corresponds to a closed curve $\gamma$ the marked point  $\gamma(0)$,
is the Serre fibration with the fibre $\Omega M$:
$$
\Lambda M \stackrel{\Omega M}{\longrightarrow} M.
$$
The exact homotopy sequence for this fibration restricted to  $\Lambda M[h]$
takes the form
\begin{equation}
\label{seq}
\dots \to \pi_i (\Lambda M [h], \omega) \stackrel{\pi_\ast}{\longrightarrow} \pi_i (M,x_0)
\stackrel{f_i}{\longrightarrow} \pi_{i-1} (\Omega_{x_0}(M), \omega) =
\end{equation}
$$
= \pi_i(M,x_0)
\to \pi_{i-1} (\Lambda M[h], \omega) \to \dots
$$
where

a) $\pi_\ast(\pi_i(\Lambda M[h], \omega)) = \mathrm{St}\,(h_i)$, where $\mathrm{St}\,(h_i)$ is the subgroup of
$\pi_i(M, x_0)$ consisting of all elements fixed under $h_i$;

b) $f_i = h_i - \mathrm{id}$ for $i \geq 2$.
}

\vskip4mm

The maps $h_i$ can be written uniformly in the simple form
$$
f_k(g) = [h,g], \ \ \ g \in  \pi_k(M,x_0), \ \ \ k \geq 1,
$$
where
$[h,g]$ is the Whitehead product of $h \in \pi_1$ and $g \in \pi_k$.

Let us consider the case when
$$
M = S^n/\Gamma, \ \ \ \ h \neq 1 \ \mbox{in $\pi_1(M, x_0)$},
$$
where $\Gamma$ acts freely and isometrically on the $n$-sphere and therefore $M^n$ is diffeomorphic to a compact space form.

If $n=2k$, then the only non-trivial group which acts freely on $S^{2k}$ is $\Z_2$ and
$S^{2k}/\Z_2 = \R P^{2k}$.

Let us consider the rational homotopy groups:
$$
\pi_i^\Q(X)  = \pi_i (X) \otimes_\Z \Q, \ \ i \geq 2.
$$
By the Cartan--Serre theorem
$$
\pi_i^\Q (S^{2k})  =
\begin{cases} \Q & \text{for $i=2k, 4k-1$} \\ 0 & \text{otherwise} \end{cases}, \ \ \
\pi_i^\Q (S^{2k+1})  =
\begin{cases} \Q & \text{for $i=2k+1$} \\ 0 & \text{otherwise} \end{cases},
$$
and moreover
$$
\pi_n (S^n) = \Z, \ \ n \geq 1.
$$

Let us assume that $n \geq 2$.

We have

\begin{theorem}
Let $M = S^n/\Gamma$, where $\Gamma$ acts freely and isometrically on $S^n$, and $h \neq 1 \in \pi_1(M)$.

Then

1)  for $i \geq 2$
$$
\pi_i^\Q (\Lambda M[h]))  =
\begin{cases} \Q & \text{for $i=4k-2, 4k-1$} \\ 0 & \text{otherwise} \end{cases}
$$
for $M = S^{2k}/\Z_2 = \R P^{2k}$
and
$$
\pi_i^\Q (\Lambda M[h]))  =
\begin{cases} \Q & \text{for $i=2k, 2k+1$} \\ 0 & \text{otherwise} \end{cases}
$$
for $M = S^{2k+1}/\Gamma$;

2)
$$
\pi_1(\Lambda M[h]) = C(h) = \Z_{r(h)} \ \ \ \mbox{for $n \geq 3$},
$$
where $C(h) =\Z_{r(h)} \subset \Gamma$ is the centralizer of $h$ in $\Gamma$,
and
$$
\pi_1 (\Lambda \R P^2 [h]) = \Z_4.
$$
\end{theorem}

{\sc Proof.}
For $h \in \pi_1(X)$ the action $h_i$ is induced by the corresponding
deck transformation of  the universal covering
$\widetilde{X} \to X$.
Therefore, if $h \neq 1$, then

\begin{enumerate}
\item
the action $h_{2k+1}$ on
$\pi_{2k+1}(S^{2k+1}/\Gamma) = \Z$ is trivial: $h_{2k+1}(z)=z$, because the deck transformation of $S^{2k+1}$
is a rotation. Therefore, by Theorem A,
$$
\pi_{2k}^\Q(\Lambda (S^{2k+1}/\Gamma)[h])  = \pi_{2k+1}^\Q(\Lambda (S^{2k+1}/\Gamma)[h])  = \Q;
$$

\item
the action $h_{2k}$ of a nontrivial element $h \in \Z_2 = \pi_1(\R P^{2k})$ on $\pi_{2k}(\R P^{2k})=\Z$ is the multiplication by $-1$:
$$
h_{2k}(z)=-z,
$$
because the corresponding deck transformation is the reflection $x \to -x$ which changes the orientation of the sphere.
It follows from Theorem A that
$$
\pi_{2k}^\Q(\Lambda \R P^{2k}[h])  = 0
$$
and
$$
\pi_{2k-1}^\Q(\Lambda\R P^{2k}[h])  = 0 \ \ \ \mbox{for $k > 1$};
$$

\item
the action $h_{4k-1}$  on
$\pi_{4k-1}(\R P^{2k})/\mathrm{Torsion}=\Z$ is trivial:
$h_{4k-1}(z)=z$,
because $\pi_{4k-1}$ is generated by the Whitehead product $[i_{2k},i_{2k}]$ where $i_{2k}$ is the generator of $\pi_{2k}$, 
$h_{2k}(i_{2k})=-i_{2k}$.
By Theorem A, we have
$$
\pi_{4k-2}^\Q(\Lambda \R P^{2k}[h])  = \pi_{4k-1}^\Q(\Lambda \R P^{2k}[h]) = \Q.
$$
\end{enumerate}

It is well-known that every commutative subgroup of $\Gamma$ is cyclic (see, for instance, \cite{Wolf}). Since therewith $\Gamma$ is finite,
then the centralizer $C(h)$ is commutative and hence cyclic: $C(h) = \Z_{r(h)}$, where $r(h)$ is the order of the maximal
cyclic subgroup of $\Gamma$ which contains $h \in \Gamma$.
By Theorem A, we have
$\pi_1(\Lambda (S^n/\Gamma)[h]) = \Z_{r(h)}$ for $n \geq 3$.

We are left to show that $\pi_1 (\Lambda \R P^2 [h]) = \Z_4$. By Theorem A, we have the exact sequence
$$
0 \to \pi_2 (\R P^2) = \Z \stackrel{\times (-2)}{\longrightarrow} \pi_2(\R P^2) = \Z \to
\pi_1 (\Lambda \R P^2[h]) \to \pi_1(\R P^2) = \Z_2 \to 0
$$
which implies the exact splitting sequence
$$
0 \to \Z_2 \to \pi_1 (\Lambda \R P^2[h],\omega) \to \Z_2 \to 0.
$$
Let us describe the homomorphisms from the last sequence. We realize  $S^2 = \{ x_1^2+ x_2^2+x_3^2=1\}$ as the unit sphere in $\R^3$. Then 
$\R P^2$ is the quotient of the unit sphere with respect to the antipodal involution.
Take the north and south poles of $S^2$: $x_3 = \pm 1$, which correspond to the same point in $\R P^2$ and
consider the paths $\kappa_\phi$ of the form
$$
x_1 = \cos \phi \sin \theta, \ x_2 = \sin\phi \sin \theta, \ x_3= \cos \theta, \ 0 \leq \theta \leq \pi.
$$
They join the poles on $S^2$ and realize the loops in $\R P^2$. This $\phi$-family of loops, where $0 \leq \phi \leq 2\pi$,
form a loop in $\Lambda \R P^2[h]$ which represents an element $[\kappa] \in \pi_1 (\Lambda \R P^2[h],\kappa_0)$.
By construction, this element generates the image of the homomorphism
$$
\pi_1(\Omega \R P^2,\kappa_0) = \pi_2(\R P^2) = \Z \to  \pi_1 (\Lambda \R P^2[h],\kappa_0)
$$
and therefore
$$
2[\kappa] = 0.
$$
Take another element $[\eta] \in \pi_1(\Omega \R P^2,\gamma_0)$
which is represented by the $\phi$-family of loops $\eta_\phi$ of the form
$$
x_1 = \sin (\theta + \phi), \ x_2 = 0, \  x_3= \cos (\theta + \phi), \ 0 \leq \theta \leq \pi,
$$
where $0 \leq \phi \leq \pi$.  It is clear that the paths $\eta_\phi$ and $\eta_{\phi+\pi}$
determine the same loop in $\R P^2$ and $\eta_0 = \kappa_0$.
By construction, the image of $[\eta]$ under the homomorphism
$$
\pi_1 (\Lambda \R P^2[h]) \to \pi_1(\R P^2) = \Z_2
$$
is nontrivial. Therefore $\pi_1 (\Lambda \R P^2[h])$ is generated by $[\eta]$ and $[\kappa]$.

Let us consider the family of paths $\eta_\phi$ with $0 \leq \phi \leq 2\pi$. It represents
$2[\eta]$ and if we take the center of each path $\eta_\phi$ and rotate the path around the axis coming through its center and the origin
we obtain an $S^1$-family of paths $\widetilde{\eta}_\phi$ which is transformed into the family $\kappa_\phi$ by a rotation of the sphere.
Hence,
$$
[\kappa] = 2[\eta]
$$
and $\pi_1(\Lambda \R P^2[h]) = \Z_4$.

This finishes the proof of Theorem 1.

\section{Homotopy quotients of the path spaces}

There is a natural $SO(2)$-action on $\Lambda M$ which consists in changing the based points.
It is described in \S \ref{sec2} where the quotient space $\Pi^+(M) = L(M)/SO(2)$ is defined.
However this action is not free, since the iterated contours have nontrivial isotropy groups.

The homotopy quotient $X_G$ of the $G$-space $X$, where $G$ is a group, is the quotient of the product
$X \times EG$ with respect to the diagonal action of $G$. Here
$$
EG \stackrel{G}{\longrightarrow} BG
$$
is the universal $G$-bundle. For $G=SO(2)$ we have $BG = \C P^\infty$, and, by definition, $EG$ is contractible.
The $G$-equivariant cohomology are defined as
$$
H^\ast_G(X) = H^\ast(X_G)
$$

Hereby we prefer to work with the spaces $LM$ which are formed by arc-length parameterised curves
and which are $SO(2)$-equivariant deformation retracts of $\Lambda M$ (see \S \ref{sec2}).
The action functional
$$
S: LM \to \R
$$
satisfies a nice property: the closed exremals of $S$ form $SO(2)$-orbits, and therefore to every nonparameterised
closed extremal of $S$ there corresponds just one critical point of the action functional
$$
S: LM_{SO(2)} \to \R
$$
(here and in the sequel, for brevity, we denote $(LM)_{SO(2)}$ by $LM_{SO(2)})$.

\begin{theorem}
Let $M = S^n/\Gamma$, where $\Gamma$ acts freely and isometrically on $S^n$, and $h \neq 1 \in \pi_1(M)$.

Then

1) for $i \geq 2$
$$
\pi_i^\Q (LM[h]_{SO(2)})  =
\begin{cases} \Q & \text{for $i=2,4k-2, 4k-1$} \\ 0 & \text{otherwise} \end{cases}
$$
for $M = S^{2k}/\Z_2 = \R P^{2k}$
and
$$
\pi_i^\Q (LM[h]_{SO(2)})  =
\begin{cases} \Q & \text{for $i=2,2k, 2k+1$} \\ 0 & \text{otherwise} \end{cases}
$$
for $M = S^{2k+1}/\Gamma$;

2)  for odd $n \geq 3$ the spaces $LM[h]_{SO(2)}$ are homotopically simple and
$$
\pi_1(LM[h]_{SO(2)}) = C(h)/\Z[h]  \ \ \ \mbox{for odd $n \geq 3$},
$$
where $C(h) \subset \Gamma$ is the centralizer of $h$ in $\Gamma$ and $\Z[h]$ is a subgroup, of $C(h)$, generated by $h$;

3) for $n \geq 1$
$$
\pi_1 (L\R P^{2n}[h]_{SO(2)}) = 0.
$$
\end{theorem}

{\sc Proof} of this theorem follows immediately from the exact homotopy sequence of the fibration
$$
LM[h] \times ESO(2) \stackrel{SO(2)}{\longrightarrow} LM[h]_{SO(2)}.
$$
It needs to clarify the statement on homotopical simplicity. This is done as follows:
for odd $n$ the fundamental groups of $S^n/\Gamma$ act trivially on the higher homotopy groups
because the deck transformations of the universal coverings are homotopy equivalent to the identity,
and now, by the explicit description of the homotopy groups of the path spaces (Theorem A),
we conclude that the path spaces and their homotopy quotients are also homotopically simple.

The groups $H_{4k-1}(L \R P^{2k}[h];\Q)$ and $H_{2k+1}(L(S^{2k+1}/\Gamma)[h];\Q)$
are generated by the following homology classes:

1) take the $(4k-1)$-dimensional manifold $V^{4k-1}$ formed by pairs $(x,v)$ where $x \in S^{2k} = \{|y|=1, y \in \R^{2k+1}\}$ and
$v$ is a unit vector tangent to $S^{2k}$ at $x$. Let us define a map $F: V^{4k-1} \to L\R P^{2k}[h]$ which
corresponds to every such a pair $(x,v)$ a semicircle $\gamma$ (in $S^{2k}$) starting at $x$ in the direction of $v$.
It is easy to show that $V^{4k-1}$ is rationally homotopy equivalent to the $(4k-1)$-sphere and
the image of the induced map in homology generates $H_{4k-1}( L\R P^{4k-1}[h];\Q)=\Q$;

2) since $h \in \Gamma$, its action on the  unit sphere $S^{2k+1}$ has the form
$$
(z_1,\dots,z_{k+1}) \to (e^{i\alpha_1}z_1,\dots,e^{i\alpha_{k+1}}z_{k+1}),
$$
where
$$
S^{2k+1} = \{z_1,\dots,z_{k+1} \in \C, |z_1|^2+\dots+|z_{k+1}|^2=1\}.
$$
Let us corresponds to every point $z \in S^{2k+1}$ the path
$$
\gamma(t) = (e^{i\alpha_1 t}z_1,\dots,e^{i\alpha_{k+1}t}z_{k+1}), \ \ 0 \leq t \leq 1,
$$
which starts at $z$ and finishes at $h(z)$.
This correspondence defines in a natural way a map
$$
F: S^{2k+1} \to L M[h]
$$
and the image of the induced map in homology generates $H_{2k+1}( LM[h];\Q)=\Q$.

Let us consider the fibrations
$$
LM[h] \times ESO(2) \stackrel{S^1}{\longrightarrow} LM[h]_{SO(2)}, \ \ \ h \neq 1 \in \pi_1(M).
$$
There are $SO(2)$-equivariant maps
$$
f: V^{4k-1} \times ESO(2) \to  L\R P^{2k}[h] \times ESO(2)
$$
and
$$
f: S^{2k+1} \times ESO(2) \to L (S^{2k+1}/\Gamma)[h] \times ESO(2)
$$
and the corresponding induced maps of the spectral sequences:
$$
f^\ast: E^n_{p,q} \to E^n_{p,q}
$$
which for $n=2$ has the forms
$$
f^\ast: H^p(L\R P^{2k}[h]_{SO(2)};H^q(S^1;\Q)) \to H^p(V^{4k-1}/SO(2);H^q(S^1;\Q)),
$$
$$
f^\ast: H^p(LM[h]_{SO(2)};H^q(S^1;\Q)) \to H^p(S^{2k+1}/SO(2);H^q(S^1;\Q)).
$$
In both cases the $SO(2)$-actions are induced by changes of the based points on the paths $\gamma$
coming into definitions of the mappings $f$. In both cases these actions are free and, in particular, we see
that

1) (for $\dim M=2k$)
$E^2_{4k-2,1} = H^{4k-2}(V^{4k-1}/SO(2);\Q)\otimes H^1(S^1;\Q)$ is generated  by $u^{k-1} \otimes v$
where $u$ is the generator of $H^2(V^{4k-1}/SO(2);\Q)$ and
$v$ is the generator of $H^1(S^1;\Q)$ such that $d_2 v = u$, 
$E^\infty_{4k-2,1} = H^{4k-1}(V^{4k-1};\Q)=\Q$ and therefore $d (u^{k-1} \otimes v) = u^{2k}=0$;

2)  (for $\dim M=2k+1$)
$E^2_{k,1} = H^{2k}(S^{2k+1}/SO(2);\Q)\otimes H^1(S^1;\Q)$ is generated  by $u^k \otimes v$
where $u$ is the generator of $H^2(S^{2k+1}/SO(2);\Q)$ and
$v$ is the generator of $H^1(S^1;\Q)$ such that $d_2 v = u$, and
$E^\infty_{2k,1} = H^{2k+1}(S^{2k+1};\Q)=\Q$ and therefore $d (u^k \otimes v) = u^{k+1}=0$.

We have
$$
u=f^\ast(u_2), \ \ \ u_2 \in H^2(LM[h]_{SO(2)};\Q),
$$
and, since
$$
E^\infty_{4k-2,1} = f^\ast (E^\infty_{4k-2,1}) = H^{4k-2}(LM[h]_{SO(2)};\Q) \otimes H^1(S^1;\Q))
$$
for $\dim M=2k$
and
$$
E^\infty_{2k,1} = f^\ast (E^\infty_{2k,1}) = H^{2k}(LM[h]_{SO(2)};\Q) \otimes H^1(S^1;\Q))
$$
for $\dim M=2k+1$,
we derive from the spectral sequences for the fibrations $LM[h] \times ESO(2) \stackrel{S^1}{\longrightarrow}
LM[h]_{SO(2)}$, that
$$
u_2^{2k} = 0 \ \ \mbox{for $\dim M=2k$ and} \ \ u_2^{k+1} = 0 \ \ \mbox{for $\dim M = 2k+1$}.
$$

Let us consider the minimal models of $LM[h]_{SO(2)}$.
The space $LM[h]_{SO(2)}$ is simply-connected for $\dim M=2k$ and is homotopically-simple with a finite cyclic
fundamental group for $\dim M=2k+1$. In both cases there are defined the minimal models
(by Sullivan; see, for instance, \cite{DGMS,FHT}). We briefly recall only the simplest properties of minimal models:

1) the minimal model $\M(X)$ of $X$ is a free graded skew-commutative algebra $\sum_{k \geq 0} \M_k$
over $\Q$ such that its generators
$\{u_\alpha\}$ (we assume that they are homogeneous: $u_\alpha \in {\mathcal M}_k$, i.e., $\deg u_\alpha=k$)
are in one-to-one correspondence with generators of $\pi^\Q(X)$;

2) there is a differential $d: \M(X) \to \M(X)$ such that
$$
d^2 = 0, \ \ d\M_k \subset \M_{k+1}
$$
and for every generator $u_\alpha$ its differential $du_\alpha$ is expressed in terms of the generators of degree less than
$\deg u_\alpha$;

3) the graded skew-commutative algebras $H^\ast(\M(X),d)$ and $H^\ast(X;\Q)$ are isomoprhic.

Let us compute the minimal models for $LM[h]_{SO(2)}$:

1) for $\dim M=2k$ the minimal model is generated by $u_2, u_{4k-2}, u_{4k-1}$ such that $\deg u_l=l, l=2,4k-2,4k-1$.
It is clear that
$$
du_2 =0, \ \ du_{4k-2} = 0,
$$
and, since as it was shown above $[u_2]^{2k}=0$ in cohomology, we have
$$
du_{4k-1} = u_2^{2k};
$$

2) for $\dim M=2k+1$ the fundamental group of $LM[h]_{SO(2)}$ is finite cyclic and therefore its rational
fundamental group is zero. The minimal model of $LM[h]_{SO(2)}$ for this nilpotent space (see \cite{DGMS,FHT})
coincide with the minimal model of its universal covering and it is generated by $u_2, u_{2k}, u_{2k+1}$ with $\deg u_l=l, l=2,2k,2k+1$.
For the same reasons as in the previous case, we have
$$
du_2 = 0, \ \ du_{2k}=0, \ \ du_{2k+1}=u_2^{k+1}.
$$

Therewith we proved the following

\begin{theorem}
Let $M = S^n/\Gamma$, where $\Gamma$ acts freely and isometrically on $S^n$, and $h \neq 1 \in \pi_1(M)$.

Then

1) for $n =2k$ the minimal model of $LM[h]_{SO(2)}$ is generated by $u_2, u_{2k}$, $u_{2k+1}$ such that
$\deg u_l=l$ for all $l$ and $du_2 = 0, du_{4k-2}=0, du_{4k-1}=u_2^{k+1}$. The cohomology ring
has the form
$$
H^\ast(LM[h]_{SO(2)};\Q) = \Q[w,z]/\{w^{2k}=0\}, \ \ \deg w=2, \deg z = 4k-2;
$$

2) for $n =2k+1$ the minimal model of $LM[h]_{SO(2)}$ is generated by $u_2, u_{2k}, u_{2k+1}$ such that
$\deg u_l=l$ for all $l$ and $du_2 = 0, du_{2k}=0, du_{2k+1}=u_2^{k+1}$. The cohomology ring
has the form
$$
H^\ast(LM[h]_{SO(2)};\Q) = \Q[w,z]/\{w^{k+1}=0\}, \ \ \deg w=2, \deg z = 2k.
$$
\end{theorem}

\section{Non-contractible closed geodesics}

Morse theory describes how closed extremals contribute to the topology of the path space.
Let us remind the result of Bott \cite{Bott}:

\medskip

{\sl Let $\gamma$ be a simple closed (Finsler or Riemann) geodesic in $M$. Then there exists a function
$$
I: S^1 = \{|z|=1, z \in \C\} \to \N
$$
such that

1) $\mathrm{ind}\, \gamma^m = \sum_{z^m=1} I(z)$, where $\mathrm{ind}$ is the Morse index of an extremal;

2) $I$ is piecewise constant and is discontinuous exactly at points $\{\lambda_1,\dots,\lambda_k\}$ which are the eigenvalues, of the complexified linearised 
Poincare mapping for $\gamma$, lying on the unit circle  $|\lambda|=1$.}

\medskip

The 2) implies that the discontinuity points $\lambda_1,\dots,\lambda_k$ are invariant with respect of the complex conjugation, and 
that there are at most $2n-2$ of them,
where $n=\dim M$ is the dimension of the (configuration) manifold.

The geodesic $\gamma^m$ is called {\it non-degenerate} if $\pm1$ does not lie in the spectrum of the corresponding Poincare mapping,
which holds if $\lambda_i^m \neq \pm 1$ for all $i=1,\dots,k$.

A metric is called {\it bumpy} if all its closed geodesics, which are different from one-point curves, are non-degenerate.

It was observed by Schwarz that an iterate $\gamma^m$ of a simple closed geodesic $\gamma$ contributes to the rational homology of $L(M)/SO(2)$
if and only if
\begin{equation}
\label{schwarz}
(\mathrm{ind}\,\gamma^m - \mathrm{ind}\,\gamma)  \ \ \mbox{is even},
\end{equation}
which, by the Bott theorem, holds if $(I(\gamma^2)-I(\gamma))$ is even or $m$ is odd.
For the equivariant cohomology we have:

\medskip

{\sl if an extremal $y$ is non-degenerate in the Morse sense, its index is equal to $k$, and $y$ corresponds either to a simple closed geodesic or
to $\gamma^m$ where $\gamma$ is a simple closed geodesic and $(\mathrm{ind}\,\gamma^m - \mathrm{ind}\,\gamma)$ is even, then
$$
H^\ast_{SO(2)}(Y^{a-\varepsilon} \cup U(y), Y^{a-\varepsilon}) = H^{\ast - k}_{SO(2)}(S^1),
$$
where $Y = LM_{SO(2)}$, $U(y)$ is a small neighborhood of $y$ in $Y$, $S(y) = a$, $\varepsilon$ is positive and sufficiently small,
and $U(y)$ contains only one critical point of $S$.}

\medskip

We refer for more details to \cite{Hingston} where the systematic application of equivariant cohomology to  closed geodesics was started
and to the survey \cite{T10} on the type numbers of closed geodesics.

In \cite{XL} the equivariant cohomology of $L \R P^{2k+1}[h]$ with coefficients in $\Z_2$ were computed and in \cite{DLX} for $\R P^3$ is was shown 
by some number-theoretical reasonings
that such cohomology can not be generated by the iterates of a single closed geodesic  of a bumpy irreversible Finsler metric.

Here we use Theorem 3 to prove analogous result for $\R P^2$:

\begin{theorem}
For every bumpy irreversible Finsler metric on $\R P^2$ there exists at least two distinct non-contractible closed geodesics.
\end{theorem}

{\sc Proof.} Let $c$ be a minimal non-contractible closed geodesic. By definition, it is simple and $\mathrm{ind}\,c = 0$. 
Its iterates $c^{2k+1}, k=0,1,\dots$, are non-contractible.
Let us assume that there are no other non-contractible closed geodesics. Then the Bott function $I(z)$  has $2 = 2n-2$ points of discontinuity, we denote them by
$e^{i\lambda}$ and $e^{-i\lambda}$.
Without loss of generality, we assume that $0 < \lambda <  \pi$. Since $\mathrm{ind}\, c=0$, we have $I(z) = 0$ for $z = e^{i\mu}$ 
with $-\lambda < \mu < \lambda$. 
By Theorem 3, 
$$
\sum_{z^{2k+1}=1} I(z) = 2l, \ \ \ k=1,2,\dots,
$$
and every even number $2l$ is presented in this form exactly twice. This implies that $I(z) = 1$ for  $z = e^{i\mu}$ with $\lambda < \mu < 2\pi - \lambda$, and, since $I(z)=0$ outside
the closure of this arc, this implies that $\lambda = \frac{\pi}{2}$ and therefore the geodesic $c^2$ is degenerate. Thus we arrive at contradiction which proves Theorem.

The following theorem demonstrates how the nontriviality of the fundamental group of $LM[h]$ can be used for proving the existence of closed geodesics.

\begin{theorem}
Let $M = S^{2n+1}/\Gamma$ and $h$ be a nontrivial element in $\pi_1(M)$. 
Let $\pi_1(LM[h])_{SO(2)} \neq 1$, $h$ has an even order in $\pi_1(M)$, and elements from $C(h)$ (the centralizer of $h$) are pairwise non-conjugate.

Then every bumpy Finsler metric on $M$ has at least two distinct closed geodesics of the class $[h]$. 
\end{theorem}

{\sc Proof.}  1) Let us assume that the metric is irreversible. 
Let $c$ be a minimal closed geodesic of $[h]$. We have $c = \gamma^k$ where $\gamma$ is a simple closed geodesic.
If all closed geodesics in $[h]$ are the iterates of $c$, then 
they are of the form $\gamma^{k+ 2pkl}$, where $l=0,1,2,\dots$, where $2p$ is the order of $h$ in $\pi_1(M)$.
By Bott theorem, the Morse indices of these iterates are even as well as $\mathrm{ind}\,c=0$.
Therefore the handle decomposition of $LM[h]_{SO(2)}$ corresponding to the action functional $S$, which is a Morse function, 
contains only even-dimensional cells. That contradicts to the non-triviality of
$\pi_1(LM[h]_{SO(2)})$ and  proves Theorem for irreversible metrics.

2) For reversible metrics, there is a possibility that $\gamma^{-m}$ is also a minimal closed geodesic for some positive $m$. However its iterates have the form 
$\gamma^{-m-2pml}, l=1,2,\dots$, and their Morse indices are also even. Hence the handle decomposition of the space still contains only even-dimensional cells 
which contradicts to the non-triviality of  $\pi_1(LM[h]_{SO(2)})$. 

Theorem is proved. 

{\sc Remark.} For reversible bumpy Finsler metrics on $\R P^{2n+1}, n=1,2,\dots$, the existence of two distinct non-contractible closed geodesics 
was established in \cite{DLX}. The following argument allows us to generalize this result for all projective spaces. 
If $c$ is a minimal geodesic, then $c^{-1}$ is also minimal. However the space $L\R P^{n}[h]/SO(2)$ is connected 
and hence there is a saddle type closed geodesic $c^\prime$ with $\mathrm{ind}\, c^\prime =1$. Since both geodesics $c$ and $c^\prime$ contribute to the
homology of the space and they have indices of different parity, by (\ref{schwarz}),  they are distinct.


\begin{thebibliography}{MMM}


\bibitem{Anosov74}
Anosov, D.V.:
Geodesics in Finsler geometry.
In: Proc. I.C.M. (Vancouver, BC 1974), vol. 2, pp. 293--297, Montreal (1975)
 (Russian. Amer. Math. Soc. Transl. {\bf 109} (1977), 81--85).
 
\bibitem{Anosov80}
Anosov, D.V.:
Some homotopies in a space of closed curves.
Math. USSR-Izv. {\bf 17} (1981), no. 3, 423--453.

\bibitem{Ballmann}
Ballmann, W.:
Geschlossene Geod\"atische auf Mannigfaltigkeiten mit unendlicher Fundamentalgruppe.
Topology {\bf 25}:1 (1986), 55--69.

\bibitem{BL}
Bangert, V., and Long, Y.:
The existence of two closed geodesics on every Finsler $n$-sphere.
Math. Ann.{\bf  346} (2010), 335--366.

\bibitem{Bott}
Bott, R.:
On the iteration of closed geodesics and the Sturm intersection theory.
Comm. Pure Appl. Math {\bf 9} (1956), 171--206.

\bibitem{DGMS}
Deligne, P., Griffiths, P., Morgan, J., and Sullivan, D.:
Real homotopy theory of K\"ahler manifolds.
Inventiones Math. {\bf 29} (1975), 245--274.

\bibitem{DLX}
Duan, H.,  Long, Y., and Xiao, Y.:
Two closed geodesics on $\R P^{2n+1}$ with a bumpy Finsler metric.
Calc. Var. Partial Diff.erential Equations {\bf 54} (2015), 2883--2894.

\bibitem{FHT}
Felix, Y., Halperin, S., and Thomas, J.-C.:
Rational Homotopy Theory, Springer-Verlag, New York, 2001.

\bibitem{Hingston}
Hingston, N.:
Equivariant Morse theory and closed geodesics.
J. Differ. Geom.  {\bf 19} (1884), 85--116.

\bibitem{Katok}
Katok, A.B.:
Ergodic perturbations of degenerate integrable Hamiltonian systems.
Math. USSR-Izvestiya {\bf 7}:3 (1973), 535--571.

 \bibitem{K}
 Klingenberg, W.:
 Lectures on closed geodesics. Springer, Berlin--Heidelberg--New York, 1978.

\bibitem{Long}
Long, Y.:
Index Theory for Symplectic Paths with Applications.
Prog. Math., vol. 207. Birkh\"auser, Basel, 2002.

\bibitem{Novikov}
Novikov, S.P.:
The Hamiltonian formalism and a many-valued analogue of Morse theory.
Russian Math. Surveys {\bf 37}:5 (1982), 1–56.

\bibitem{T86}
Taimanov, I.A.:
Closed geodesics on non-simply-connected manifolds.
Russian Math. Surveys {\bf 40}:6 (1986), 143--144.

\bibitem{T0}
Taimanov, I.A.:
Closed extremals on two-dimensional manifolds.
Russian Math. Surveys {\bf 47}:2 (1992), 163--211.

\bibitem{T10}
Taimanov, I.A.:
The type numbers of closed geodesics.
Regular and Chaotic Dynamics {\bf 15} (2010), 84--100.

\bibitem{Wolf}
Wolf, J.A.:
Spaces of Constant Curvature. University of California, Berkeley, California, 1972.

\bibitem{XL}
Xiao, Y., and Long, Y.:
Topological structure of non-contractible loop space and closed geodesics on
real projective spaces with odd dimensions.
Adv. Math. {\bf 279} (2015), 159--200.



\end{thebibliography}
\end{document}